\documentclass[10pt]{amsart}
\usepackage{amsthm,mathrsfs,graphicx,mathabx,tcolorbox,color,tikz,enumitem,hyperref,cancel,graphics,bm}
\usepackage[normalem]{ulem}
\usepackage{amsmath,bbm,tabu}
\usepackage{amssymb,verbatim}
\usepackage{url}

\definecolor{orange}{rgb}{0.898, 0.621, 0.0}
\definecolor{skyblue}{rgb}{0.336, 0.703, 0.910}
\definecolor{bluishgreen}{rgb}{0, 0.617, 0.449}
\definecolor{yellow}{rgb}{0.937, 0.890, 0.258}
\definecolor{blue}{rgb}{0, 0.445, 0.695}
\definecolor{red}{rgb}{0.832, 0.367, 0}
\definecolor{purple}{rgb}{0.797, 0.473, 0.652}

\newtheorem{theorem}{Theorem}
\numberwithin{theorem}{section}
\newtheorem{proposition}[theorem]{Proposition}

\newtheorem{conjecture}[theorem]{Conjecture}
\newtheorem{question}[theorem]{Question}
\theoremstyle{definition}
\newtheorem{definition}[theorem]{Definition}

\begin{document}

\begin{abstract}
In this paper, we present a minimal counterexample to a conjecture of Perles that answers a question of Haase and Ziegler. The example is a simple 4-polytope that has an induced 3-connected 3-regular subgraph, whose graph complement is connected. This subgraph is planar and not the graph of a facet of the polytope.
\end{abstract}

\title{A minimal counterexample to a strengthening of Perles' Conjecture}
\author{Joseph Doolittle}

\maketitle

\section{Introduction}

The question of reconstructibility of a polytope from its vertex-edge graph has been of great interest to many. One of the earliest attempts at making progress on this problem was a conjecture by Perles, that the facets of simple polytopes were precisely the induced, 3-connected, non-separating, \((d-1)\)-regular subgraphs of the vertex-edge graph \cite{Pe70}. If the conjecture were true, then simple polytopes could be reconstructed from their vertex-edge graphs. Blind and Mani, Kalai, and Friedman showed that all simple polytopes can be reconstructed from their graph \cite{BM87, Ka88, Fr07}. However, their proofs did not settle Perles' Conjecture. The conjecture was latter shown to be false in general by Haase and Ziegler \cite{HZ02}, although their counterexample leaves something to be desired. The authors themselves note that the construction is ``as-concrete-as-possible.'' Furthermore, in this counterexample, the \(3\)-regular subgraph of interest in their example is not planar, as all graphs of facets of 4-polytopes are.

The example of Haase and Ziegler was explicitly realized by Witte \cite{Wi02}. Its f-vector is (2592, 5194, 3121, 529). We present two new counterexamples to Perles' Conjecture, summarized here by their f-vectors. The minimum f-vector of 4-polytopes with a planar counterexample to Perles' conjecture is (38, 76, 50, 12). The minimum f-vector of simple 3-spheres with a counterexample to a generalization of Perles' conjecture is (38, 76, 49, 11). Minimality is considered first with regard to number of facets, then with regard to number of vertices. We only consider minimality within the class of either 3-spheres or 4-polytopes, as appropriate, and the examples presented realize each of these minima.

Section \ref{defs} gives definitions for the terms used. A minimal planar counterexample is given concretely in Section \ref{example}. In Section \ref{spheres}, we give more information about the verification of minimality, as well as the extension to simplicial spheres. In Section \ref{story} we describe some of the constructions that led to this example. In Section \ref{compute}, we describe which parts of this paper were done through computation and the code used.

\section{Definitions}\label{defs}

While most of these definitions are standard, we reiterate them here for completeness.

\begin{definition}
A \emph{d-polytope} is a subset of \(\mathbb{R}^d\) that is the convex hull of a finite set of points in \(\mathbb{R}^d\). A \emph{d-polytope} is a compact subset of \(\mathbb{R}^d\) that is the intersection of a finite set of half-spaces of \(\mathbb{R}^d\). That these definitions are equivalent is a celebrated and important result in polytope theory.
\end{definition}

\begin{definition}
An \emph{\(i\)-face} of a polytope is a set of points in a polytope that maximize a linear functional and whose affine span is \(i\)-dimensional. A \((d-1)\)-face of a \(d\)-polytope is called a facet, and a \((d-2)\)-face of a \(d\)-polytope is called a ridge.
\end{definition}

\begin{definition}
The \emph{polar dual} of a polytope with the origin in its interior is the polytope defined by the intersection of the half-spaces \(\{\overrightarrow{y} : \overrightarrow{y}\cdot\overrightarrow{x} \leq 1\}\), where \(\overrightarrow{x}\) are the coordinates of the 0-faces of the polytope.
\end{definition}

\begin{definition}
A \emph{simplicial complex} is a collection of sets \(\Delta\), such that if \(\sigma \in \Delta\) and \(\tau \subset \sigma\), then \(\tau \in \Delta\). The elements of \(\Delta\) of size \(i+1\) are called \emph{\(i\)-faces}. We represent simplicial complexes by only the maximal faces, which we also call \emph{facets}. When all the facets are the same size, we call the complex \emph{pure}. When the complex is pure, then \emph{ridges} are the faces of cardinality one less than the facets.
\end{definition}

\begin{definition}
The \emph{star} of a face \(\sigma\) in a simplicial complex \(\Delta\) is the simplicial complex generated by the facets containing \(\sigma\), \(st_\Delta(\sigma)= \{\tau \mid \tau \subset F\text{ for some } F \in \Delta \text{ such that } \sigma \subset F \} \). The \emph{link} of a face \(\sigma\) in a simplicial complex is the set of faces whose intersection with \(\sigma\) is empty, and whose union with \(\sigma\) is in \(\Delta\), \(lk_\Delta(\sigma) = \{\tau \mid \tau \cap \sigma = \emptyset, \tau \cup \sigma \in \Delta\}\).
\end{definition}

\begin{definition}
A \emph{vertex-edge graph} of either a simplicial complex or a polytope is the graph whose vertices are the 0-dimensional faces and whose edges are the pairs of 0-dimensional faces which are contained in a 1-dimensional face.
\end{definition}

\begin{definition}
A \emph{simple} \(d\)-polytope is one such that its vertex-edge graph is \(d\)-regular.
\end{definition}

\begin{definition}
A \emph{facet-ridge graph} of either a simplicial complex or a polytope is the graph whose vertices are the facets and whose edges are pairs of facets that contain a common ridge. This definition is not the only possible definition for facet-ridge graph, but it is sufficient for our purposes and nicely reflects the vertex-edge definition.
\end{definition}

\begin{definition}
The \emph{boundary} of a \(d\)-dimensional simplicial complex is the simplicial complex generated by the faces of dimension less than \(d\) that are contained in at most one \(d\)-dimensional face.
\end{definition}

\begin{definition}
An \emph{interior facet} of a simplicial complex is a facet for which none of its ridges are in the boundary of the simplicial complex.
\end{definition}

\begin{definition}
A \emph{psuedomanifold} with boundary is a topological space which has a triangulation such that the triangulation is pure, every ridge is in at most two facets, and the facet-ridge graph is connected.
\end{definition}

\begin{definition}
A \emph{pinched \(i\)-face} in a pure simplicial complex of dimension \(d\) is an \(i\)-face such that the geometric realization of its link is homeomorphic to \(\mathbb{S}^{d-i-1} \times [0,1]^i\).
\end{definition}

This definition is intentionally more broad than what will be used in this paper. For our purposes, we only consider pinched vertices in a 3-dimensional complex. Then the condition to be a pinched vertex is that its link is \(\mathbb{S}^1 \times [0,1]\), or an annulus. For example, the complex generated by the facets \(0123,0234,0345,0456,0561,0612\) has \(0\) as a pinched vertex, since its link is an annulus. The motivation for the name is that this example is like taking a piece of dough and pinching it with two fingers in the middle.

\begin{definition}
A \emph{Perles piece} is a connected simplicial complex such that every facet contains exactly one ridge in the boundary of the complex.
\end{definition}

\begin{conjecture}[Perles]
Every \((d-1)\)-regular, induced, \((d-1)\)-connected, non-separating subgraph of the vertex-edge graph of a simple \(d\)-polytope corresponds to a facet of \(P\).
\end{conjecture}

The boundary of the dual of a simple polytope is a simplicial complex. Under duality, facets become vertices, and vertices become facets. Consequently, the subgraph in the conjecture above becomes a subgraph of the facet-ridge graph under duality. That the subgraph is \((d-1)\)-regular implies that every facet of the dual has a ridge in common with \(d-1\) other facets, and has one ridge that is not in any other facet, which is therefore in the boundary. That the subgraph is induced is equivalent to downward containment, which is a requirement of simplicial complexes by their definition. We can then restate Perles' Conjecture.

\begin{conjecture}[Perles]
The stars of vertices are the only non-separating Perles pieces in the polar dual of a simple polytope.
\end{conjecture}

At the end of \cite{HZ02}, which gave the first counterexample to Perles' Conjecture, the authors ask the following question.

\begin{question}[Haase and Ziegler]
Are the stars of vertices the only non-separating Perles pieces with planar facet-ridge graphs in the polar dual of a simple 4-polytope?
\end{question}

We answer this question in the next section.

\section{A minimal planar counterexample to Perles' Conjecture}\label{example}

In this section, we give an explicit description of a facet-minimal example of a simple 4-polytope whose graph has a planar, \(3\)-regular, induced, \(3\)-connected, non-separating subgraph that does not correspond to a facet. For ease of writing and consistency, we only refer to the polar dual of the example.
We describe the example by the following list of twelve vertices, which correspond to the twelve facets of the facet-minimal example.

\[\begin{tabu}{rrrrrl}
A=( &  27, & - 95, &  120, &    0 & )     \\
B=( & -50, & - 45, &  101, & - 94 & )  \\
C=( & - 9, & - 67, &  126, & - 35 & )    \\
D=( & 195, & -145, &   11, &  125 & )  \\
E=( & -40, & - 10, &    8, & - 65 & )    \\
F=( & 232, & -102, & - 21, &  198 & ) \\
G=( & -63, & - 25, &   94, & -139 & )  \\
H=( & -80, &   45, & - 49, & - 65 & )   \\
I=( & -72, &    4, &   24, & - 90 & )     \\
J=( & -30, &  167, & -154, &   92 & )  \\
K=( & -43, &  190, & -199, &  100 & ) \\
L=( & -67, &   80, & - 61, & - 26 & )    \\
\end{tabu}\]
\[P = \text{conv}(\{A,B,\ldots,L\})\]

The polytope \(P\) is simplicial, so its boundary forms a simplicial complex. The following is the facet-ridge graph of a non-separating Perles piece contained in this complex.

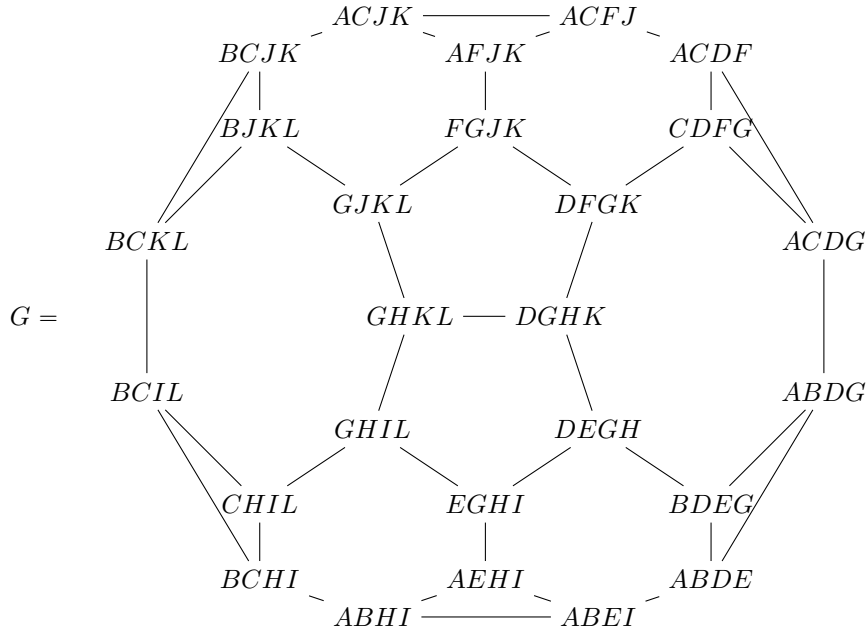
\begin{tikzpicture}
\node at (-6  , 0  ) {$G=$};
\node at (-1  , 0  ) (v1){$GHKL$};
\node at ( 1  , 0  ) (v2){$DGHK$};
\node at ( 0  , 2.5) (v3){$FGJK$};
\node at ( 0  ,-2.5) (v4){$EGHI$};
\node at (-1.5, 1.5) (v5){$GJKL$};
\node at ( 1.5, 1.5) (v6){$DFGK$};
\node at (-1.5,-1.5) (v7){$GHIL$};
\node at ( 1.5,-1.5) (v8){$DEGH$};
\node at ( 3  , 2.5) (v9){$CDFG$};
\node at ( 3  , 3.5) (v10){$ACDF$};
\node at ( 1.5, 4  ) (v11){$ACFJ$};
\node at ( 0  , 3.5) (v12){$AFJK$};
\node at (-1.5, 4  ) (v13){$ACJK$};
\node at (-3  , 3.5) (v14){$BCJK$};
\node at (-3  , 2.5) (v15){$BJKL$};
\node at ( 3  ,-2.5) (v16){$BDEG$};
\node at ( 3  ,-3.5) (v17){$ABDE$};
\node at ( 1.5,-4  ) (v18){$ABEI$};
\node at ( 0  ,-3.5) (v19){$AEHI$};
\node at (-1.5,-4  ) (v20){$ABHI$};
\node at (-3  ,-3.5) (v21){$BCHI$};
\node at (-3  ,-2.5) (v22){$CHIL$};
\node at ( 4.5, 1  ) (v23){$ACDG$};
\node at (-4.5, 1  ) (v24){$BCKL$};
\node at ( 4.5,-1  ) (v25){$ABDG$};
\node at (-4.5,-1  ) (v26){$BCIL$};

\draw (v1) -- (v2) -- (v6) -- (v3) -- (v5) -- (v1) -- (v7) -- (v4) -- (v8) -- (v2);
\draw (v6) -- (v9) -- (v10) -- (v11) -- (v12) -- (v13) -- (v14) -- (v15) -- (v5);
\draw (v12) -- (v3);
\draw (v13) -- (v11);
\draw (v7) -- (v22) -- (v21) -- (v20) -- (v19) -- (v18) -- (v17) -- (v16) -- (v8);
\draw (v19) -- (v4);
\draw (v20) -- (v18);
\draw (v9) -- (v23) -- (v25) -- (v16);
\draw (v15) -- (v24) -- (v26) -- (v22);
\draw (v14) -- (v24);
\draw (v26) -- (v21);
\draw (v10) -- (v23);
\draw (v17) -- (v25);
\end{tikzpicture}

This planar graph is not the graph of the star of a vertex, which corresponds to a subgraph which is not the graph of a facet in the dual. This answers the question of Haase and Ziegler in the negative. This example is a minimal such example with respect to the number of vertices of \(P\), which was verified by exhaustive computer search of all simplicial spheres with fewer vertices.

\section{Simplicial Spheres}\label{spheres}
In the previous section, we gave a counterexample to Perles' Conjecture where the 3-connected graph is planar. In this section, we relax the planar condition and expand the class of objects we consider from 4-polytopes to simplicial 3-spheres. These relaxations allow a counterexample with one fewer vertex, from 12 vertices in the dual of the previous example, to 11 vertices in this example.

Kalai posed the following conjecture: All simplicial spheres can be reconstructed from their facet-ridge graph \cite{Ka09}. For polytopal simplicial spheres, under duality this is equivalent to the result of Blind and Mani, Kalai, and Friedman. The barrier to extending these results is that convexity can no longer be assumed, since each proof at some point requires some geometric information that is derived from convexity. Although Perles' Conjecture is false even for polytopes, it has a natural dual extension to simplicial spheres: ``The stars of vertices are the only Perles pieces of a simplicial sphere.''
 
Here, we provide a counterexample to such a conjecture, which is minimal with respect to the number of vertices, then minimal with respect to the number of facets. We give the facet list of a vertex-minimal simplicial sphere that contains a non-separating Perles piece that is not the star of a vertex.

The facets of a Perles piece that is not the star of a vertex are

\begin{tabular}{lllllll}
ABCD & BCDF & BFGI & CEGI & CEIK & DEFH & DGIK \\
ABDE & BCFG & BGHI & CEHJ & CFGI & DFHJ & DGIJ \\
ACDF & BDEH & BGHJ & CEJK & CHJK & DFJK & DGJK \\
ADEF & BEHJ & EGIK & FHJK & GHIJ \\
\end{tabular}

Adding the following facets to those above forms a simplicial sphere.

\begin{tabular}{lllllll}
ABCE & BCEG & BDHI & CEFH & CFIK \\
ACEF & BDFI & BEGJ & CFHK & EGJK \\
DFIK & DHIJ \\
\end{tabular}

This sphere is shellable, since it is the boundary of the following shellable 4-ball, whose facets are given in a shelling order, left to right, then top to bottom.

\begin{tabular}{lllllll}
ADFJK & ADFIK & ABDFI & AFHJK & ACFHK & ACHJK & ACEJK \\
ADFHJ & ACFIK & ACFGI & ACEIK & ACEGI & AEGIK & ADGIK \\
ADGJK & ADGIJ & ADHIJ & AGHIJ & ABDHI & ABGHI & ABGHJ \\
AEGJK & ABEGJ & ABEHJ & ABCEG & ACEHJ & ACEFH & ADEFH \\
ABFGI & ABCFG & ABCDF & ABDEH \\
\end{tabular}

The facet-ridge graph of the Perles piece has the following subgraph:

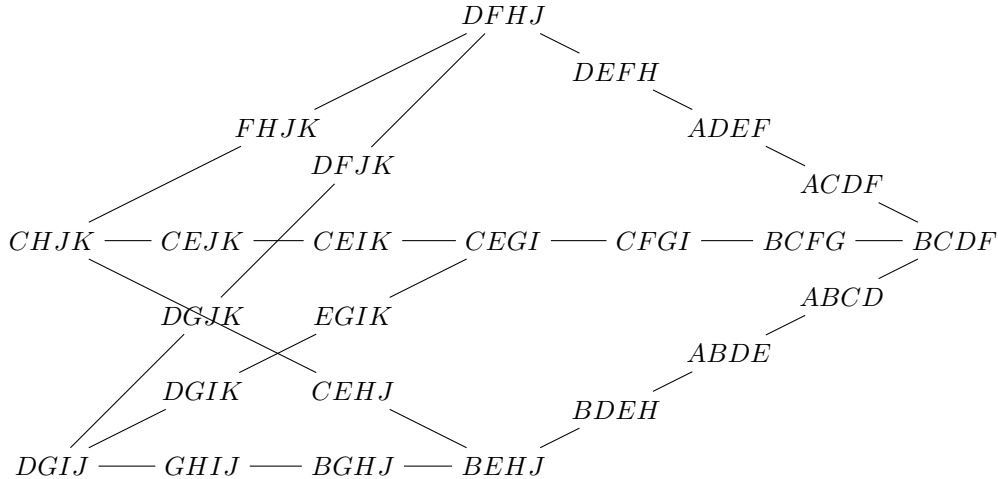
\begin{tikzpicture}
\node at (-3  ,-3  ) (v1){$DGIJ$};
\node at (-3  , 0) (v2){$CHJK$};
\node at ( 9  , 0) (v3){$BCDF$};

\node at ( 3  ,-3  ) (v4){$BEHJ$};
\node at ( 3  , 0) (v5){$CEGI$};
\node at ( 3  , 3) (v6){$DFHJ$};

\node at (-1  ,-3) (v7){$GHIJ$};
\node at ( 1  ,-3) (v8){$BGHJ$};
\node at (-1  ,-2) (v9){$DGIK$};
\node at ( 1  ,-1) (v10){$EGIK$};
\node at (-1  ,-1) (v11){$DGJK$};
\node at ( 1  , 1) (v12){$DFJK$};

\node at ( 1  , -2) (v13){$CEHJ$};
\node at (-1  , 0) (v14){$CEJK$};
\node at ( 1  , 0) (v15){$CEIK$};
\node at ( 0  , 1.5) (v16){$FHJK$};

\node at ( 7.5, -.75) (v17){$ABCD$};
\node at ( 6  ,-1.5 ) (v18){$ABDE$};
\node at ( 4.5,-2.25) (v19){$BDEH$};
\node at ( 7  , 0  ) (v20){$BCFG$};
\node at ( 5  , 0  ) (v21){$CFGI$};
\node at ( 7.5, 0.75) (v22){$ACDF$};
\node at ( 6  , 1.5 ) (v23){$ADEF$};
\node at ( 4.5, 2.25) (v24){$DEFH$};

\draw (v1) -- (v7) -- (v8) -- (v4);
\draw (v1) -- (v9) -- (v10) -- (v5);
\draw (v1) -- (v11) -- (v12) -- (v6);

\draw (v2) -- (v13) -- (v4);
\draw (v2) -- (v14) -- (v15) -- (v5);
\draw (v2) -- (v16) -- (v6);

\draw (v3) -- (v17) -- (v18) -- (v19) -- (v4);
\draw (v3) -- (v20) -- (v21) -- (v5);
\draw (v3) -- (v22) -- (v23) -- (v24) -- (v6);

\end{tikzpicture}

This graph has \(K_{3,3}\) as a minor, and therefore the Perles piece does not have a planar facet-ridge graph.

The facet-ridge graph of the Perles piece that is not the star of a vertex in this example is not planar. This example is not unique, and all examples with 11 vertices are available by email upon request.
It is not currently known if any of these spheres are polytopal. These examples were generated by computer search, which further verified that all simplicial 3-spheres with 10 or fewer vertices do not have counterexamples to the extension of Perles' Conjecture to simplicial spheres.

\begin{question}
Are any of the 11-vertex spheres with a Perles piece that is not the star of a vertex polytopal?
\end{question}

\section{Methods}\label{story}

Bohus, Jockusch, Lee and Prabhu showed that any Perles piece not the star of a vertex must not be a topological ball or topological torus \cite{BJLP98}. In the same paper, they remark that Rudin's ball comes close to being a Perles piece, although they incorrectly identify how many interior facets it has. Of the 41 facets of Rudin's ball, 17 are interior \cite{Ru58}. A relatively more tractable triangulation of the ball with similar properties is Ziegler's ball \cite{Zi98}.  In this triangulation, only 5 of the facets are interior, even closer to a Perles piece.

We now describe a pair of general operations that combine to form a useful modification to simplicial complexes. We require that in all cases these operations are done to subcomplexes which satisfy the conditions of psuedomanifolds when considered as a part of the larger complex.

The first operation is simply stellar subdivision of a facet. This takes a general facet \(1234\) and replaces it with the facets \(1235, 1245, 1345, 2345\).

The second operation requires a collection of facets isomorphic to \(1234, 1245, 1256\) with \(123\) a boundary ridge, and both of \(1245\) and \(1256\) interior facets. We replace these facets with the facets \(1734, 7234, 1745, 7245, 1256\). After this operation all of these facets have a boundary triangle, respectively \(173, 723, 175, 725, 125\).

We combine these operations in the following way. We find a collection of facets isomorphic to \(0123, 1234, 2345\) where the facet isomorphic to \(1234\) is an interior facet, and the ridges isomorphic to \(012\) and \(345\) are in the boundary. We then stellar subdivide the facet isomorphic to \(1234\) with a vertex we will label \(6\). Then we perform the second operation under the isomorphism that takes the facets \(0123, 1236, 1246\) to the facets \(1234, 1245, 1256\). We perform the second operation a second time, this time under the isomorphism that takes the facets \(1634, 6234, 2345\) to the facets \(1234, 1245, 1256\). This combination of three operations replaces three facets with ten facets, and each of these ten facets have a boundary ridge. This also creates a pinched vertex at \(6\).

It turns out that the 5 interior facets of Ziegler's ball all satisfy the conditions to apply this operation. After applying the operation 5 times, we obtain a simplicial complex each of whose facets have exactly one boundary triangle. This complex has 57 facets and 25 vertices. We consider how this complex is different from a topological ball, which can never have this property \cite{BJLP98}. The topological difference is the 5 pinched vertices, the vertices used in subdividing the interior facets.

We may then focus on constructions that use these pinched vertices as a means to the end, rather than a consequence of the construction. The following facet list describes a basic building block that we will use in the following constructions.

\(\Gamma = \langle 0125, 0256, 0236, 0346, 0134, 0145, 1345, 2346 \rangle\). 

If we were to attach additional facets along \(246, 256, 125, 135\), then each of the facets of \(\Gamma\) have exactly one boundary triangle. We create three copies of \(\Gamma\) with vertices \(0_a,1_a,2_a,3_a,4_a,5_a,6_a\), \(0_b,1_b,2_b,3_b,4_b,5_b,6_b\), \(0_c,1_c,2_c,3_c,4_c,5_c,6_c\). We then glue these copies together by making the identifications \(1_a=2_b\), \(3_a=4_b\), \(5_a=6_b\), \(1_b=2_c\), \(3_b=4_c\), \(5_b=6_c\), \(1_c=2_a\), \(3_c=4_a\), \(5_c=6_a\). Then we cone all the faces that are copies of the original \(256\) and \(125\) faces to a new vertex, \(7\). In this way we have constructed a complex, so that each of the facets in the copies of \(\Gamma\) have exactly one boundary triangle, and the 6 facets containing \(7\) also each have only one boundary triangle. Overall, this example uses 4 pinched vertices, \(0_a,0_b,0_c,7\), 13 vertices and 30 facets, substantially smaller than our previous example.

We may do the same construction with three copies of \(\Gamma\), but then do some additional subdivisions and identifications to allow the triangles we previously coned to be glued together instead. This leads to the example in Section \ref{example}, which only has 3 pinched vertices, 12 vertices and 28 facets.

\begin{question}
Are there any Perles pieces with 1 or 2 pinched vertices?
\end{question}

Each of these gives an \(\mathbb{R}^3\)-embeddable simplicial complex. We may then add more facets to complete these to a sphere. In the case that this is also the facet list of a simplicial polytope, we may take its dual to get a polytope with a counterexample to Perles' Conjecture.

\section{Computation}\label{compute}

There are several claims made in this paper that rely on computation. We will go through claim by claim and detail how the claim was verified.

\begin{proposition}
No simplicial 3-sphere with 10 vertices contains a Perles piece that is not the star of a vertex.
\end{proposition}

\begin{proposition}
There are 1669 simplicial 3-spheres with 11 vertices that contain a Perles piece that is not the star of a vertex. 
\end{proposition}

These two propositions were proved with the same code. The computations for the first proposition were done on the author's home computer, an Intel Core 2 Q9400 processor. The computations for the second proposition were done on the KU High Performance Computing cluster.

The code for these computations is available by email upon request. It consists of three main parts. The first part is lextet version 0.24 \cite{lextet}. This generates all 3-manifolds with a given number of vertices. The second part and third part are both written by the author. The second part computes the rank of homology of the manifold, and makes sure that the ranks of homology matches that of the sphere. The last part finds all Perles pieces of the simplicial complex, using a depth-first search algorithm. Facets are included or excluded from the candidate Perles piece, and the branch terminates when it finds that a facet of the candidate Perles piece is forced to have 0, 2, 3 or 4 boundary ridges.

Separately, Sage was used to check homology of the resulting manifolds, and found some non-sphere manifolds, which were thrown out \cite{sage}. The remaining manifolds total 1669 and are available by email upon request.

\begin{proposition}
The Perles pieces of those 1669 simplicial spheres that are not the stars of vertices do not have planar facet-ridge graphs.
\end{proposition}

Using just the third part of the code described above, the non-star Perles pieces were identified for each manifold. Sage was used to check planarity of each of the facet ridge graphs, and they are all not planar \cite{sage}.

\begin{proposition}
The coordinates given in Section \ref{example} describe a polytope containing the specified Perles piece.
\end{proposition}

While this proposition can be independently verified with any number of software packages, constructing coordinates is a difficult problem. The author recommends polymake for verification of the example \cite{polymake}. The code which found these coordinates does so in the following way. First, it computes the  determinants of the matrices of the chirotope derived from the simplicial sphere. It then adjusts the position of the vertices based on those determinants to attempt to make more of them positive. When all of the determinants are positive, the coordinates describe a polytope whose boundary is the desired simplicial sphere. This is a closely related process to spring algorithms to generate nice embeddings of graphs. This is non-deterministic, and is helpful only for providing coordinates for a polytopal realization of a simplicial sphere. It will not prove that a simplicial sphere is not polytopal. The initial input is generated by finding an edge within the sphere whose link is large, then choosing coordinates that naturally induce all the facets containing that edge. For more information on realizing simplicial spheres as polytopes and chirotopes, the author recommends a paper by Firsching \cite{Fi17}.

\section{Acknowledgments}\label{thanks}

The author would like to thank their advisor Marge Bayer for her continued support. The author would like to thank the High Performance Computing cluster at KU and Weizhang Huang for help getting the code completed in a reasonable amount of time. The author would like to thank Eran Nevo for conversations and inspiration on these and related problems. The author would like to thank Micheal Joswig for catching an error in the coordinates in Section \ref{example}.

\bibliographystyle{siam}
\bibliography{GlobalBib}

\end{document}